\documentclass{article}
\usepackage{amssymb,amsmath}
\usepackage{graphicx}

\def\qed{\hfill {\hbox{${\vcenter{\vbox{               
   \hrule height 0.4pt\hbox{\vrule width 0.4pt height 6pt
   \kern5pt\vrule width 0.4pt}\hrule height 0.4pt}}}$}}}

\def\tr{\triangleright}
\def\bar{\overline}

\newtheorem{theorem}{Theorem}
\newtheorem{definition}{Definition}

\newtheorem{proposition}[theorem]{Proposition}
\newtheorem{corollary}[theorem]{Corollary}
\newtheorem{example}{Example}

\newenvironment{proof}[1][Proof]{\smallskip\noindent{\bf #1.}\quad}%
{\qed\par\medskip}

 \date{}

\title{\Large \textbf{On bilinear biquandles}}

\author{\begin{tabular}{c} Sam Nelson \\ 
\small Department of Mathematics \\
\small  Pomona College \\
\small 610 North College Avenue \\ 
\small Claremont, CA 91711 \\
\texttt{knots@esotericka.org}\end{tabular}
\and
\begin{tabular}{c}  Jacquelyn L. Rische \\ 
\small Department of Mathematics \\
\small 103 MSTB \\
\small University of California, Irvine \\
\small Irvine, CA 92697 \\
\texttt{jrische@math.uci.edu}\end{tabular}
}

\begin{document}

\maketitle

\begin{abstract}
We define a type of biquandle which is a generalization of symplectic 
quandles. We use the extra structure of these bilinear biquandles to define
new knot and link invariants and give some examples.
\end{abstract}



\section{\large \textbf{Introduction}}

A \textit{biquandle} is an algebraic structure consisting of a set $B$ with 
four binary operations
$(a,b)\mapsto a^b, a_b, a^{\bar{b}}$ and $a_{\bar{b}}$ satisfying
axioms derived from the oriented Reidemeister moves, where generators
of the algebra are identified with semi-arcs in an oriented link diagram. 
A biquandle may also be understood as a solution $S:B\times B \to B\times B$ 
to the \textit{set-theoretic Yang-Baxter equation} 
\[(S\times \mathrm{Id})(\mathrm{Id}\times S)(S\times \mathrm{Id}) =
(\mathrm{Id}\times S)(S\times \mathrm{Id})(\mathrm{Id}\times S)\]
which satisfies some additional criteria corresponding to the first and
second Reidemeister moves. Such a map $S$ is called a \textit{switch}, and
a biquandle is an invertible switch with components $S(a,b)=(b_a,a^b)$ 
satisfying $S^{-1}(a,b)=(b^{\bar{a}},a_{\bar{b}})$ and the extra conditions 
required by the reverse type II and type I moves. 

\begin{figure}
\[\scalebox{1.4}{\includegraphics{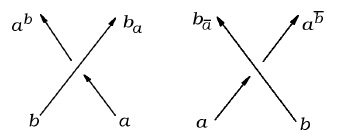}}\]
\caption{Biquandle operations at crossings}
\end{figure}

This relationship between the biquandle axioms and the Reidemeister moves
makes biquandles a natural source of knot and link invariants. For example,
the biquandle counting invariant $|\mathrm{Hom}(B(L),T)|$ is the cardinality
of the set of biquandle homomorphisms from the knot biquandle of a link
$L$ into a finite target biquandle $T$. One can think of each homomorphism
as a ``coloring'' of the link diagram by $T$, assigning an element of $T$ to
every semiarc in a diagram of $L$ such that the biquandle operations are
satisfied at every crossing; we can then see that this family of invariants 
is a generalization of Fox's $n$-coloring invariants. Indeed, biquandles 
generalize quandles which in turn generalize knot groups. Biquandles 
have been studied in recent papers such as \cite{FJK}, \cite{K1}, \cite{C1} 
and more.

Many of the examples of biquandles in the current literature are natural 
generalizations of types 
of quandle structures -- Alexander biquandles generalize Alexander quandles, 
Silver-Williams switches generalize Joyce's homogeneous quandles, etc. 
In this paper we generalize the symplectic quandles studied in \cite{NN}
(also known as \textit{quandles of transvections}) to define what we call
bilinear biquandles.

The paper is organized as follows. In section 2 we list the biquandle axioms
and give examples of biquandles. In section 3 we define bilinear biquandles, 
obtain some results about their structure and give an example of bilinear 
biquandle which is not a quandle. In section 4 we generalize the symplectic 
quandle polynomial invariants defined in \cite{NN} to the biquandle case and 
give some examples of classical and virtual links which have the same value
for the biquandle counting invariant but are distinguished by the bilinear
biquandle invariant. In section 5, we list all bilinear biquandle structures 
with cardinality up to 27 as determined by our computer search. In section 6, 
we end with some questions for further research.

\section{\large \textbf{Biquandles and symplectic quandles}} \label{s2}

Let $B$ be a set. A \textit{biquandle structure} on $B$ consists of four 
binary operations $(a,b)\mapsto a^b, a_b, a^{\bar{b}}$ and $a_{\bar{b}}$ 
such that

\begin{list}{}{}
\item[1.]{For every $a,b\in B$ we have
\[a=a^{b\bar{b_a}}, \quad b=b_{a\bar{a^b}}, \quad a=a^{\bar{b}b_{\bar{a}}},
\quad \mathrm{and}
\quad b=b_{\bar{a}a^{\bar{b}}},\]}
\item[2.]{for every $a,b\in B$ there exist $x,y\in B$ such that
\[ x=a^{b_{\bar{x}}},\quad  a=x^{\bar{b}},\quad  b=b_{\bar{x}a},\quad
y=a^{\bar{b_y}},\quad a=y^b,\quad  \mathrm{and} \quad b=b_{y\bar{a}},\]}
\item[3]{for every $a,b,c\in B$ we have
\[a^{bc}=a^{c_bb^c}, \quad c_{ba}=c_{a_bb_a}, \quad 
(b_a)^{c_{a^b}}=(b^c)_{a^{c_b}},\]
\[a^{\bar{b}\bar{c}}=a^{\bar{c_{\bar{b}}}\bar{b^{\bar{c}}}}, \quad 
c_{\bar{b}\bar{a}}=c_{\bar{a_{\bar{b}}}\bar{b_{\bar{a}}}}, \quad 
\mathrm{and}\quad 
(b_{\bar{a}})^{\bar{c_{\bar{a^{\bar{b}}}}}}
=(b^{\bar{c}})_{\bar{a^{\bar{c_{\bar{b}}}}}},\]}
\item[4.]{for every $a\in B$ there exist $x,y\in B$ such that
\[x=a_x, \quad a=x^a, y=a^{\bar{y}}\quad \mathrm{and} \quad a=y_{\bar{a}}.\]}
\end{list}

These axioms are obtained from the oriented Reidemeister moves by thinking
of each semiarc (portion of the knot diagram between over/under crossing 
points) as a biquandle element; the biquandle elements on the outside of
each pictured diagram portion must then agree before and after the move. For 
example, the oriented type III move with all positive crossings is below:

\[\scalebox{2.3}{\includegraphics{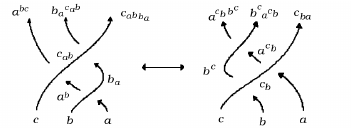}}\]

While there are eight oriented type III moves, in the presence of the 
oriented type II moves we need only the two type III moves with all positive 
and all negative crossings. See \cite{K1} for more.

\begin{example}
\textup{A \textit{quandle} is a set $Q$ with binary operations 
$\tr,\tr^{-1}:Q\times Q\to Q$ such that 
\begin{list}{}{}
\item[(i)]{for all $a\in Q$, $a\tr a=a,$}
\item[(ii)]{for all $a,b\in Q$ we have $(a\tr b)\tr^{-1}b=a=(a\tr^{-1}b)\tr b$
and}
\item[(iii)]{for all $a,b,c\in Q$ we have $(a\tr b)\tr c=(a\tr c)\tr 
(b\tr c).$}
\end{list}
Every quandle is a biquandle with $a^b=a\tr b$, $a^{\bar{b}}=a\tr^{-1} b,$ 
and $a_b=a_{\bar{b}}=a.$
}
\end{example}

\begin{example}
\textup{As an example of a non-quandle biquandle, let $B=\mathbb{Z}_n$ and
let $s,t\in \mathbb{Z}_n$ be invertible elements. Then $B$ is a biquandle 
with}
\[a^b=ta+(1-st)b, \quad a^{\bar{b}}=t^{-1}a+(1-s^{-1}t^{-1})b, \quad a_b=sa,
\quad \mathrm{and} \quad a_{\bar{b}}=s^{-1}a.\]
\end{example}

\begin{example}
\textup{More generally, let $M$ be any module over $\mathrm{Z}[s^{\pm 1},
t^{\pm 1}]$.
Then $M$ is a biquandle under the operations defined in example 2, called an
\textit{Alexander biquandle}. See \cite{K1} and \cite{LN} for more.}
\end{example}

\begin{example}
\textup{Let $D$ be an oriented link diagram, i.e. a planar 4-valent graph with
two inward-oriented and two outward-oriented edges incident on every vertex, 
with vertices decorated to indicate crossing information. Then the 
\textit{knot biquandle} $B(L)$ of the link represented by $L$ has a 
presentation with one generator for each edge and relations at each crossing 
as depicted  in figure 1. The elements of the knot biquandle are equivalence 
classes of biquandle words in these generators under the equivalence relation
generated by the biquandle axioms and the crossing relations.
For example, the trefoil knot below has the listed 
knot biquandle.}
\[\raisebox{-0.5 in}{\includegraphics{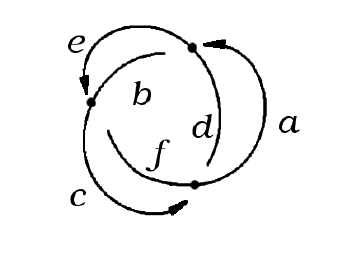}} \quad
\begin{array}{rcl} 
\langle a,b,c,d,e,f \ & | \ &  a^d=b, b_e=c, c^f=d, \\
 & & d_a=e, e^b=f, f_c=a \rangle.\end{array}
\]
\textup{If we drop the planarity requirement, such a $D$ defines a 
\textit{virtual link}, and we obtain a knot biquandle by the same procedure.
Crossings arising from non-planarity are depicted as circled intersections;
see \cite{KV} for more about virtual knots and links. }
\end{example}

If $B=\{x_1,\dots, x_n\}$ is a finite biquandle, we can represent $B$ 
symbolically with a $2n\times 2n$ block matrix, where the four $n\times n$ 
blocks encode the four biquandle operations. That is, we let
\[M_B=\left[
\begin{array}{c|c} 
M^1 & M^2 \\ \hline
M^3 & M^4
\end{array}
\right]
\quad M^l_{ij} = k \quad \mathrm{where} \quad x_k=\left\{
\begin{array}{ll}
(x_i)^{\bar{(x_j)}} & l=1 \\
(x_i)^{(x_j)} & l=2 \\
(x_i)_{\bar{(x_j)}} & l=3 \\
(x_i)_{(x_j)} & l=4 \\
\end{array}
\right.
\]

For example, the Alexander biquandle $B=\mathbb{Z}_3$ with $s=2,t=1$ has 
biquandle matrix
\[M_B=\left[
\begin{array}{ccc|ccc}
3 & 2 & 1 & 3 & 2 & 1 \\
1 & 3 & 2 & 1 & 3 & 2 \\
2 & 1 & 3 & 2 & 1 & 3 \\ \hline
2 & 2 & 2 & 2 & 2 & 2 \\
1 & 1 & 1 & 1 & 1 & 1 \\
3 & 3 & 3 & 3 & 3 & 3 \\
\end{array}
\right].\]

Biquandle matrices can be used to compute the biquandle counting invariant
$|\mathrm{Hom}(B(K),T)|$ and the Yang-Baxter 2-cocycle invariants of a
knot or link in an algebra-agnostic way; see \cite{NV} and \cite{CN} for more.

\begin{example}\textup{
As a final example of a non-quandle biquandle structure, let $R$ be a 
ring with identity and $M$ an $R$-module. Then the operations
\[y_x=Ax+By,\quad  x^y=Cx+Dy \]
with $A,B\in R$ and $C=A^{-1}B^{-1}A(1-A)$ and $D=1-A^{-1}B^{-1}AB$
define an \textit{invertible switch} $S:M\times M \to M\times M$ by
$S(x,y)=(y_x,x^y)$ provided $[B,(A-I)(A,B)]=0$ where
$[X,Y]=XY-YX$ and $(X,Y)=X^{-1}Y^{-1}XY$; such a switch gives a 
biquandle structure with barred operations 
$S^{-1}(x,y)=(x^{\bar{y}},y_{\bar{x}})$ provided the axioms arising from 
the type I and reverse type II moves are satisfied. See \cite{BF} for more.}
\end{example}

\section{\large \textbf{Bilinear biquandles}} \label{res}

Let $R$ be any commutative ring and $M$ any free module over $R$. Let 
$\langle,\rangle:M\times M\to R$ be an antisymmetric bilinear form.
Then $M$ is a biquandle with operations 
\[\mathbf{x}^{\bar{\mathbf{y}}}=\mathbf{x}-
\langle \mathbf{x},\mathbf{y}\rangle \mathbf{y}, \quad
\mathbf{x}^{\mathbf{y}}=\mathbf{x}+
\langle \mathbf{x},\mathbf{y}\rangle \mathbf{y}, \quad
\mathbf{x}_{\bar{\mathbf{y}}}=\mathbf{x} \quad \mathrm{and} \quad
\mathbf{x}_{\mathbf{y}}=\mathbf{x}.
\]
Such a biquandle is in fact a quandle; quandles of this type have been called
\textit{symplectic quandles} or \textit{quandles of transvections}. See 
\cite{NN} and \cite{Z}.

We would like to extend this definition to define non-quandle biquandles. 

\begin{definition}
\textup{Let $M$ be a free module over a commutative ring $R$. A} bilinear 
biquandle 
\textup{structure on $M$ is a biquandle structure on $M$ such that}
\[\mathbf{x}^{\mathbf{y}}
=\alpha \mathbf{x}+f(\mathbf{x},\mathbf{y})\mathbf{y}, \quad
\mathbf{x}^{\bar{\mathbf{y}}}
=\alpha' \mathbf{x}+f'(\mathbf{x},\mathbf{y})\mathbf{y}, \quad
 \mathbf{x}_{\mathbf{y}}=\beta \mathbf{x}
\quad \mathrm{and} \quad
\mathbf{x}_{\bar{\mathbf{y}}}=\beta'\mathbf{x},
\]
\textup{where $\alpha,\alpha',\beta,\beta'\in R$ and $f,'f:M\times M\to R$
are bilinear forms.}
\end{definition} 

We know already that $\alpha=\alpha'=\beta=\beta'=1$, $f'(x,y)=-f(x,y)$,
$f$ an antisymmetric bilinear form gives us a biquandle structure, 
namely the symplectic quandle structure just described. We would like to know,
then, what other bilinear biquandles are possible?

We start with some observations.

\begin{proposition}
Let $M$ be a bilinear biquandle. Then $\alpha'=\alpha^{-1}$ and 
$\beta'=\beta^{-1}.$
\end{proposition}

\begin{proof}
Here we consider the biquandle axioms arising from the direct type II move:
\[\scalebox{1.4}{\includegraphics{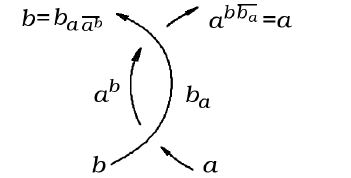}}\]
Since $b_a=\beta b$, we have $b_{a\bar{a^b}}=\beta'(\beta b)=b$, and thus
$\beta\beta'=1$.

Moreover, since $a^{b\bar{b_a}}=a$ for all $a,b\in M$, taking $b=0\in M$
we have $a^0=\alpha a$, $0_a=0$ and $a^{\bar{0}}=\alpha'a$, and thus 
$a^{0\bar{0_a}}=\alpha \alpha'a=a$ and $\alpha\alpha'=1$.
\end{proof}

\begin{proposition}
For any bilinear biquandle, we must have $f(a,a)=\beta^{-1}-\alpha$ and
$f'(a,a)=\beta-\alpha^{-1}$.
\end{proposition}

\begin{proof}
Consider the biquandle axiom derived from the Reidemeister type I move
with a positive crossing:
\[\includegraphics{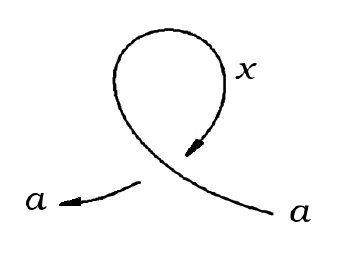} \raisebox{0.5in}{$a_x=x, \quad x^a=a$}\]
Here $a_x=\beta a=x$, and we note that this $x$ is unique since $\beta$ is
invertible; then 
\[x^a=(\beta a)^a
= \alpha (\beta a)+f(\beta a, a) a 
= (\alpha\beta  +\beta f(a,a))a=a\]
and hence $1=\alpha\beta +\beta f(a,a)$. Then 
\[\beta f(a,a)=1-\alpha\beta \Rightarrow f(a,a)=\beta^{-1}-\alpha.\]

The other case is similar.
\end{proof}

\begin{corollary}
Let $A\in M_m(R)$ be the matrix of $f$ with respect to an ordered basis 
$\{\mathbf{b_1},\dots,\mathbf{b_m}\}$ of $B$, so that 
\[f(\mathbf{x},\mathbf{y})=\mathbf{x}A\mathbf{y}^t \quad 
\mathrm{where} \quad \displaystyle{\mathbf{x}=
\sum_{i=1}^m x_i\mathbf{b_i}}.\]
Then the diagonal entries of $A$ must satisfy $A_{ii}=\beta^{-1}-\alpha$. 
\end{corollary}

\begin{proposition}
In any bilinear biquandle, we must have 
\[f'(\mathbf{x},\mathbf{y})=\omega f(\mathbf{x},\mathbf{y}) \quad 
\mathrm{where} \quad \omega
=-\alpha^{-2}\beta^{-2}-\alpha^{-1}\beta+\alpha^{-2}.\]
\end{proposition}

\begin{proof}
Using the direct type II move, we have 
$\mathbf{x}=\mathbf{x}^{\mathbf{y}\bar{\mathbf{y}_{\mathbf{x}}}}$. Therefore,
\begin{eqnarray*} 
\mathbf{x} 
& = & 
(\alpha \mathbf{x}+f(\mathbf{x},\mathbf{y})\mathbf{y})^{\bar{\beta \mathbf{y}}}
\\
& = & \alpha^\prime(\alpha \mathbf{x}+f(\mathbf{x},\mathbf{y})\mathbf{y})
     +f^\prime(\alpha \mathbf{x}+f(\mathbf{x},\mathbf{y})\mathbf{y},\beta 
\mathbf{y})\beta \mathbf{y} \\
& = & \alpha^\prime \alpha \mathbf{x} + \alpha^\prime f(\mathbf{x},
\mathbf{y})\mathbf{y}
+f^\prime(\alpha \mathbf{x}, \beta \mathbf{y})\beta \mathbf{y}+
f^\prime(f(\mathbf{x},\mathbf{y})\mathbf{y},\beta \mathbf{y})\beta \mathbf{y}.
\end{eqnarray*}
We know that $\alpha^\prime=\alpha^{-1},$ so with further simplification, 
we get 
\[\mathbf{x}=\mathbf{x}+\alpha^{-1}f(\mathbf{x},\mathbf{y})\mathbf{y}+
\alpha \beta^2f^\prime(\mathbf{x},\mathbf{y})\mathbf{y}
+\beta^2f(\mathbf{x},\mathbf{y})f^\prime(\mathbf{y},\mathbf{y})\mathbf{y}.\] 
We have $f^\prime(\mathbf{y},\mathbf{y})=\beta-\alpha^{-1}$, so 
\[\mathbf{x}=\mathbf{x}+\alpha^{-1}f(\mathbf{x},\mathbf{y})\mathbf{y}
+\alpha \beta^2f^\prime(\mathbf{x},\mathbf{y})\mathbf{y}
+\beta^2f(\mathbf{x},\mathbf{y})(\beta-\alpha^{-1})\mathbf{y}.\] 
Then 
\[0=\alpha^{-1}f(\mathbf{x},\mathbf{y})\mathbf{y}
+\alpha \beta^2f^\prime(\mathbf{x},\mathbf{y})\mathbf{y}
+\beta^2f(\mathbf{x},\mathbf{y})(\beta-\alpha^{-1})\mathbf{y},\] 
and since this is true for all $\mathbf{y}\in B$, we must have 
\[0=\alpha^{-1}f(\mathbf{x},\mathbf{y})+\alpha \beta^2f^\prime(\mathbf{x},
\mathbf{y})+\beta^2f(\mathbf{x},\mathbf{y})(\beta-\alpha^{-1}).\] 
Hence, 
\[-\alpha \beta^2f'(\mathbf{x},\mathbf{y})=(\alpha^{-1}+
\beta^2(\beta-\alpha^{-1}))f(\mathbf{x},\mathbf{y}).\] 
So, 
\begin{eqnarray*} 
f'(\mathbf{x},\mathbf{y})
& = & (-\alpha^{-1}\beta^{-2})(\alpha^{-1}+\beta^{3}-
\beta^2\alpha^{-1})f(\mathbf{x},\mathbf{y}) \\
& = & (-\alpha^{-2}\beta^{-2}-\alpha^{-1}\beta+
\alpha^{-2})f(\mathbf{x},\mathbf{y}).
\end{eqnarray*}
\end{proof}

The type III Reidemeister move axioms impose conditions on the bilinear form
$f(\mathbf{x},\mathbf{y})=\mathbf{x}A\mathbf{y}^t$:

\begin{proposition} \label{r3}
In any bilinear biquandle $B$, we must have
\[
\alpha(1-\beta^2)f(\mathbf{x},\mathbf{y}) =0 \quad \mathrm{and}  \quad
\beta(1-\beta^2)f(\mathbf{x},\mathbf{y}) =0.
\]
In particular, the entries $A_{ij}$ of the matrix $A$ such that 
$f(\mathbf{x},\mathbf{y})=\mathbf{x}A\mathbf{y}^t$ with respect to 
a basis of $B$ must satisfy 
\[\alpha(1-\beta^2)A_{ij}=\beta(1-\beta^2)A_{ij}=0\]
for all $i,j.$
\end{proposition}

\begin{proof}
The middle strand in the general case gives us the equation 
$(b_a)^{c_{a^b}}={b^c}_{a^{c_b}}$. Then
\[(b_a)^{c_{a^b}} = (\beta b)^{(\beta c)} 
= \alpha (\beta b) + f(\beta b, \beta c) \beta c 
= \alpha \beta b + \beta^3 f(b,c) c\]
while
\[(b^c)_{a^{c_b}}  =  \beta (\alpha b + f(b,c) c) 
= \alpha\beta b +\beta f(b,c) c
\]
so we must have $\beta f(b,c)=\beta^3 f(b,c)$ for all $b,c\in B.$

Now, we note that 
$\mathbf{x}^0=\alpha \mathbf{x}+f(\mathbf{x},0)0=\alpha \mathbf{x}$.
The undercrossing strand gives us the equation $a^{bc}=a^{b_cc^b}$, so the 
special case $c=0$ says
\[a^{b0}=\alpha(a^b)=\alpha(\alpha a+f(a,b) b)=\alpha^2 a+\alpha f(a,b)b\]
while
\[a^{b_00^b}=a^{(\beta b) 0}=\alpha(\alpha a+f(a,\beta b)\beta b)=
\alpha^2 a +\alpha\beta^2f(a,b)b\]
and hence $\alpha f(a,b)=\alpha\beta^2f(a,b)$ for all $a,b\in B.$
\end{proof}

These observations constrain the possible bilinear biquandle structures
on $(\mathbb{Z}_n)^m$ enough to make it practical to find all such biquandle 
structures for small values of $n$ and $m$ by computer search. Specifically,
given invertible $\alpha,\beta\in\mathbb{Z}_n$, we compute the corresponding 
list of all $x\in\mathbb{Z}_n$ satisfying $\alpha(1-\beta^2)x=
\beta(1-\beta^2)x=0$; these are candidates for entries in an $m\times m$
matrix $A$, with diagonal entries $\beta^{-1}-\alpha$. We then compute the 
biquandle operation matrix for each triple $(\alpha,\beta,A)$ over the set
$(\mathbb{Z}_n)^m$ and test the resulting operation
matrix for the biquandle axioms, rejecting any triples which fail to 
satisfy all of the axioms. \texttt{Maple} code
implementing this procedure for $(\mathbb{Z}_n)^2$ and $(\mathbb{Z}_n)^3$ 
is available in the file \texttt{bilinear-biquandles.txt} downloadable from 
\texttt{www.esotericka.org}. The results for $(\mathbb{Z}_2)^2$, 
$(\mathbb{Z}_3)^2$, $(\mathbb{Z}_4)^2$, $(\mathbb{Z}_5)^2$, 
$(\mathbb{Z}_2)^3$, $(\mathbb{Z}_3)^3$ and $(\mathbb{Z}_2)^4$  
are collected in section \ref{lst}.

\begin{example} \label{bb1}
\textup{Let $T=(\mathbb{Z}_4)^2$ and let $\alpha=\beta=3$, 
$\omega=-\alpha^{-2}\beta^{-2}-\alpha^{-1}\beta+\alpha^{-2}=1$ and 
$A=\left[\begin{array}{cc} 0 & 2 \\ 2 & 0 \end{array}\right]$. Then
one checks that the operations
\[\mathbf{x}^{\mathbf{y}}=
\mathbf{x}^{\bar{\mathbf{y}}}=
3\mathbf{x}+\mathbf{x}
\left[\begin{array}{cc} 0 & 2 \\ 2 & 0 \end{array}\right]
\mathbf{y}^t\mathbf{y}, \quad
\quad \mathbf{x}_{\mathbf{y}}=3\mathbf{x} = \mathbf{x}_{\bar{\mathbf{y}}}
\] define a bilinear biquandle structure on $T$.}

\textup{For example, let $\mathbf{x}=(x_1,x_2)$ and $\mathbf{y}=(y_1,y_2)$.
Then after a bit of arithmetic we find that}
\begin{eqnarray*}
\mathbf{x}^{\mathbf{y}\bar{\mathbf{y}_{\mathbf{x}}}}
& = & (x_1 +20x_2y_1^2+20y_1x_1y_2+72y_2x_2y_1^3+72y_1^2x_1y_2^2, \\
& & x_2+20y_2x_2y_1+20x_1y_2^2+72y_2^2x_2y_1^1+72y_1x_1y_2^3)
\end{eqnarray*}
\textup{which is just $(x_1,x_2)$ in $(\mathbb{Z}_4)^2$.}
\end{example}

\section{\large \textbf{Link invariants from bilinear biquandles}}

Since a bilinear biquandle is not just a biquandle but also an $R$-module,
we can take advantage of this extra structure to enhance the biquandle
counting invariant as in \cite{NN}.

\begin{definition}\textup{
Let $L$ be a link, $B(L)$ the knot biquandle of $L$ and $T$ a finite bilinear 
biquandle. Define the} bilinear biquandle polynomial \textup{ of $L$ with
respect to $T$ to be}
\[\phi_{BB}(L,T,q,z)=\sum_{f\in\mathrm{Hom}(B(L),T)} 
q^{|\mathrm{Im}(f)|}z^{|\mathrm{Span}(\mathrm{Im}(f))|}.\]
\end{definition}

Since $\phi_{BB}$ is determined by the set $\mathrm{Hom}(B(L),M)$ which is
an invariant of link type, so is $\phi_{BB}$. In particular, $\phi_{BB}$ 
specializes to the biquandle counting invariant $|\mathrm{Hom}(B(L),M)|$
when $q=z=1$, though in general $\phi_{BB}$ contains more information
than the counting invariant alone.

\begin{example}
\textup{ Let $T$ be the bilinear biquandle defined in example \ref{bb1}.
Then the pictured virtual 
link $L$ is distinguished from the trefoil knot $3_1$ by the associated
bilinear biquandle invariant $\phi_{BB}$, though both links have the
same counting invariant value.}
\[\includegraphics{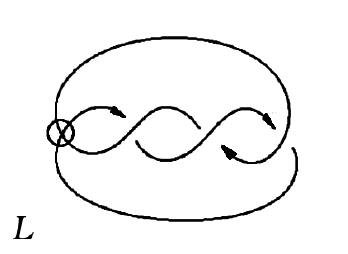} 
\raisebox{0.5in}{$
\begin{array}{c}
\phi_{BB}(L,T,q,z)=qz+6q^2z^2+3qz^2+6q^2z^4 \\ \\
|\mathrm{Hom}(B(L),T)|=16\end{array}$}\]
\[\includegraphics{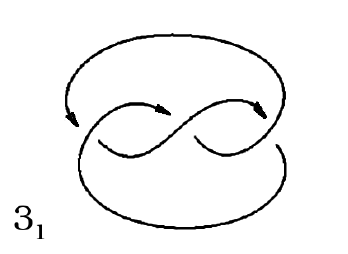} 
\raisebox{0.5in}{$
\begin{array}{c}\phi_{BB}(3_1,T,z,q)=qz+3qz^2+12q^2z^4 \\ \\
|\mathrm{Hom}(B(L),T)|=16\end{array}$}\]
\textup{Note that if we specialize $z=1$, ignoring the module structure and
using only the biquandle structure, the resulting invariant fails to 
distinguish the virtual links as both reduce to $q+3q+12q^2$.}
\end{example}

\begin{example}\textup{
Each monomial in a value of $\phi_{BB}$ corresponds to a sub-biquandle of
the target biquandle, namely the image of the knot biquandle under some
homomorphism. If two knots or links have different values of $\phi_{BB}$,
the difference in these values can yield information about the difference
between the knots not apparent from the counting invariant alone. Let
$T$ be the bilinear biquandle $T=(\mathbb{Z}_4)^2$ with}
\[\mathbf{x}^{\mathbf{y}}=\mathbf{x}^{\bar{\mathbf{y}}}=
\mathbf{x}+
\mathbf{x}\left[\begin{array}{cc} 2 & 1 \\ 1 & 2 \end{array}\right]
\mathbf{y}^t\mathbf{y},\quad \quad \quad
\mathbf{x}_{\mathbf{y}} =\mathbf{x}_{\bar{\mathbf{y}}} =3\mathbf{x}
\]
\textup{Then the links $L_1$ and $L_2$ have invariant values
$\phi_{BB}(L_1)=qz+48q^3z^4+6q^2z^2+30q^2z^4+72q^4z^8+3qz^2$ and
$\phi_{BB}(L_2)=qz+6q^2z^2+30q^2z^4+3qz^2$ respectively. Thus, there are
sub-biquandles of $B_1,B_2\subset T$ with cardinality $3$ and $4$ respectively
such that $L_1$ has biquandle colorings by $B_1$ and $B_2$ while $L_2$ does
not. Moreover, the submodule spanned by $B_1$ has cardinality 4, while
the submodule spanned by $B_2$ has cardinality 8.}
\[\includegraphics{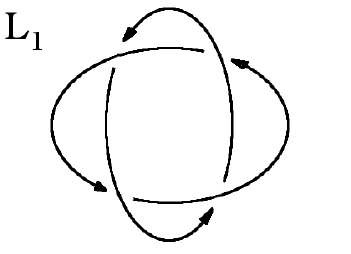} \quad 
\includegraphics{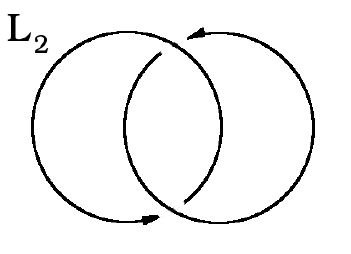}\]
\end{example}

\section{\large \textbf{Bilinear biquandles of small cardinality}} \label{lst}

In table 1 we list the results of our computer search for bilinear 
biquandles of small cardinality. These results were obtained using 
\texttt{Maple} programs in the file \texttt{bilinear-biquandles.txt}
available at \texttt{www.esotericka.org/quandles}. In light of the results
of section \ref{res}, we identify each bilinear biquandle by listing
$\alpha,\ \beta$ and the matrix $A$ of $f(\mathbf{x},\mathbf{y})$ with
respect to the standard basis of $(\mathbb{Z}_n)^m$. We list only those
bilinear biquandles which are not symplectic quandles and which have
cardinality less than or equal to 27.

\begin{table}[!h]
\begin{center}
\begin{tabular}{|cccc|cccc|cccc|} \hline
 & & & & & & & & & & & \\ 
$(\mathbb{Z}_n)^m$ & $\alpha$ & $\beta$ & $A$ & 
$(\mathbb{Z}_n)^m$ & $\alpha$ & $\beta$ & $A$ & 
$(\mathbb{Z}_n)^m$ & $\alpha$ & $\beta$ & $A$ \\ 
 & & & & & & & & & & & \\ \hline
 & & & & & & & & & & & \\ 
$(\mathbb{Z}_3)^2$ & 2 & 2 & 
$\left[\begin{array}{cc} 0 & 0 \\ 0 & 0 \end{array}\right]$ & 
$(\mathbb{Z}_3)^2$ & 2 & 2 & 
$\left[\begin{array}{cc} 0 & 1 \\ 2 & 0 \end{array}\right]$ & 
$(\mathbb{Z}_4)^2$ & 1 & 3 & 
$\left[\begin{array}{cc} 2 & 0 \\ 2 & 2 \end{array}\right]$ \\
& & & & & & & & & & & \\ 
$(\mathbb{Z}_4)^2$ & 1 & 3 & 
$\left[\begin{array}{cc} 2 & 1 \\ 1 & 2 \end{array}\right]$ &
$(\mathbb{Z}_4)^2$ & 3 & 1 & 
$\left[\begin{array}{cc} 2 & 0 \\ 2 & 2 \end{array}\right]$ &
$(\mathbb{Z}_4)^2$ & 3 & 1 & 
$\left[\begin{array}{cc} 2 & 1 \\ 1 & 2 \end{array}\right]$ \\
& & & & & & & & & & & \\ 
$(\mathbb{Z}_4)^2$ & 3  & 3 & 
$\left[\begin{array}{cc} 0  & 0 \\ 0 & 0 \end{array}\right]$ &
$(\mathbb{Z}_4)^2$ & 3 & 3 & 
$\left[\begin{array}{cc} 0  & 2 \\ 2 & 0 \end{array}\right]$ &
$(\mathbb{Z}_4)^2$ & 3 & 3 & 
$\left[\begin{array}{cc} 0  & 1 \\ 3 & 0 \end{array}\right]$ \\
& & & & & & & & & & & \\ 
$(\mathbb{Z}_5)^2$ & 4 & 4 & 
$\left[\begin{array}{cc} 0 & 0 \\ 0 & 0 \end{array}\right]$ & 
$(\mathbb{Z}_5)^2$ & 4 & 4 & 
$\left[\begin{array}{cc} 0 & 1 \\ 4 & 0 \end{array}\right]$ &
$(\mathbb{Z}_3)^3$ & 2  & 2 & 
$\left[\begin{array}{ccc} 
0 & 0 & 0 \\ 
0 & 0 & 0 \\ 
0 & 0 & 0 \end{array}\right]$ \\
& & & & & & & & & & & \\ 
\hline
\end{tabular}
\end{center}
\caption{Non-quandle bilinear biquandles of cardinality $\le 27$}
\end{table}

\section{\large \textbf{Questions}}

In this section we collect a few questions for future research.

Our initial computer search did not turn up any examples of finite 
biquandles in which all four operations have the form 
\[(\mathbf{x},\mathbf{y})\mapsto \alpha_i \mathbf{x} + 
f_i(\mathbf{x},\mathbf{y})\mathbf{y}\]
where $f_i:B\times B\to R$ is a nonzero bilinear form for all $i=1,2,3,4$. 
Are there any examples of such biquandles?

We notice from table 1 that for $n=2,3,5$ the only bilinear biquandle 
structures on $(\mathbb{Z}_n)^m$ we seem to find are only slight variations 
of the symplectic quandle structure -- the bilinear form is antisymmetric and 
$\alpha,\beta\in \{-1,1\}$. Does this pattern hold for all prime $n$?

There are, of course other possible combinations of module elements and 
bilinear forms similar to the symplectic quandle structure that one could
use in looking for finite biquandles, such as
\[(\mathbf{x},\mathbf{y})\mapsto f_i(\mathbf{x},\mathbf{y})\mathbf{x} +
g_i(\mathbf{x},\mathbf{y})\mathbf{y},\quad i=1,2,3,4\]
or 
\[(\mathbf{x},\mathbf{y})\mapsto A_i\mathbf{x} +
g_i(\mathbf{x},\mathbf{y})\mathbf{y},\quad i=1,2,3,4\]
where $f_i,g_i$ are bilinear forms and $A_i\in M_m(R)$. All of these 
should have $\phi_{BB}$ type invariants associated. Which of these 
formats give interesting new finite biquandles?

\end{document}